\documentclass[12pt]{article}
\usepackage{amsmath, amsthm, amssymb}
\numberwithin{equation}{section}
\begin{document}
\author{Ajai Choudhry}
\title{Ideal solutions of the Tarry-Escott Problem\\ of degree seven}
\date{}
\maketitle

\begin{abstract} 
In this paper we obtain four new parametric ideal solutions of the Tarry-Escott problem of degree 7, that is, of the simultaneous diophantine equations, $
\sum_{i=1}^8x_i^r=\sum_{i=1}^8y_i^r,\;r=1,\,2,\,\dots,\,7$. While all the known parametric solutions of the problem, with one exception, are given by polynomials of degrees $ \geq 5$, the solutions obtained in this paper are given by quartic polynomials, and are thus simpler than almost all of the known solutions.

\end{abstract}

\smallskip
  

\section{Introduction}\label{Intro}

\normalsize
The Tarry-Escott problem (written briefly as TEP) of degree $k$ consists of finding two distinct multisets
of integers $\{x_1,\,x_2,\,\ldots,\,x_s\}$ and $\{y_1,\,y_2,\,\ldots,\,y_s\}$ such that
\begin{equation}
\sum_{i=1}^sx_i^r=\sum_{i=1}^sy_i^r,\;\;\;r=1,\,2,\,\dots,\,k,
\label{tepsk}
\end{equation}
where $k$ is a given positive integer. When $s=k+1$, solutions of the diophantine system \eqref{tepsk} are known as ideal solutions. Parametric ideal solutions of the TEP are known only when $k \leq 7$. 

This paper is concerned with  parametric ideal solutions the Tarry-Escott problem of degree 7, that is, of the diophantine system,
\begin{equation}
\sum_{i=1}^8x_i^r=\sum_{i=1}^8y_i^r,\;\;\;r=1,\,2,\,\dots,\,7.
\label{tep7}
\end{equation}
 Except for a single  parametric ideal solution that is expressible by  quadratic polynomials (see \cite[p. 633]{Che}), all other known parametric ideal solutions of \eqref{tep7} are given by univariate polynomials of degrees $\geq 5$.  Recently MacLeod \cite{ML1, ML2} has examined various  methods of solving the diophantine system \eqref{tep7}  and applied the theory of elliptic curves to obtain several parametric ideal solutions of degree 5. An alternative method give by Choudhry  \cite[pp.\ 410--412]{Cho2} also yields a parametric ideal solution of degree 5. The complete ideal solution of the TEP of degree 7 is not yet known.

In this paper we give an elementary  method of obtaining ideal solutions of the TEP of degree 7. We apply the method to obtain four parametric ideal solutions of degree 4. This is an improvement over the known results since, as mentioned above, all known parametric solutions, with one exception, are of degrees $> 4$.  We also  show how more parametric ideal solutions of degrees $> 4$ may be obtained.

\section{Parametric ideal solutions of the TEP of degree 7}\label{parmsol}
To solve the diophantine system \eqref{tep7}, we write
\begin{equation}
 x_{i+4}=-x_i,\;\; y_{i+4}=-y_i,\;\;i=1,\,\dots, 4, \label{simplecond}
\end{equation}
when \eqref{tep7} is identically satisfied for odd values of the exponent $r$, and the diophantine system reduces to 
\begin{equation}
\sum_{i=1}^4x_i^r=\sum_{i=1}^4y_i^4,\;\;\;r=2, 4, 6.
\label{tep7red}
\end{equation} 
In fact, all known solutions of the TEP of degree 7 are obtained by applying the simplifying conditions \eqref{simplecond}.

We note that Choudhry \cite[p.\ 119]{Cho1} has given the following solution of the diophantine system $\sum_{i=1}^4x_i^r=\sum_{i=1}^4y_i^4,\;\;\;r=1, 2, 4$:
\begin{equation}
\begin{aligned}
x_1 &= \alpha_1m + (\alpha_1 + 2\alpha_3)(\alpha_1 + 2\alpha_2 + \alpha_3)n, \\
x_2 &= \alpha_2m - (\alpha_1 - \alpha_3)(2\alpha_1 + \alpha_2 + 2\alpha_3)n, \\
x_3 &= \alpha_3m - (2\alpha_1 + \alpha_3)(\alpha_1 + 2\alpha_2 + \alpha_3)n,\\
 x_4 &=- (\alpha_1 + \alpha_2 + \alpha_3)m + (\alpha_1 - \alpha_3)(\alpha_1 - \alpha_2 + \alpha_3)n, \\
y_1 &= -\alpha_1m + (\alpha_1 + 2\alpha_3)(\alpha_1 + 2\alpha_2 + \alpha_3)n, \\
y_2 &= -\alpha_2m - (\alpha_1 - \alpha_3)(2\alpha_1 + \alpha_2 + 2\alpha_3)n, \\
y_3 &= -\alpha_3m - (2\alpha_1 + \alpha_3)(\alpha_1 + 2\alpha_2 + \alpha_3)n, \\
y_4 &= (\alpha_1 + \alpha_2 + \alpha_3)m + (\alpha_1 - \alpha_3)(\alpha_1 - \alpha_2 + \alpha_3)n,
\end{aligned}
\label{RMJMsol}
\end{equation}
where $\alpha_1,  \alpha_2, \alpha_3, m$ and $n$ are arbitrary parameters. Accordingly, we choose $x_i, y_i, i=1,\ldots, 4$, as given by \eqref{RMJMsol} when \eqref{tep7red} is identically satisfied for $r=2$ and $r=4$, while for $r=6$, we may write \eqref{tep7red}, after transposing all terms to one side, as follows:
\begin{multline}
12mn(\alpha_1^2 - \alpha_3^2)(\alpha_1 + \alpha_2)(\alpha_1 + 2\alpha_2 + \alpha_3)(\alpha_2 + \alpha_3)\{m^2-9(\alpha_1 + \alpha_3)^2n^2\}\\
\times \{(2\alpha_1^2 - \alpha_1\alpha_2 + 5\alpha_1\alpha_3 - \alpha_2^2 - \alpha_2\alpha_3 + 2\alpha_3^2)m^2 - (2\alpha_1^4 - \alpha_1^3\alpha_2 \\
+ 9\alpha_1^3\alpha_3 - \alpha_1^2\alpha_2^2+ 37\alpha_1^2\alpha_2\alpha_3 + 14\alpha_1^2\alpha_3^2 + 38\alpha_1\alpha_2^2\alpha_3 \\
+ 37\alpha_1\alpha_2\alpha_3^2 + 9\alpha_1\alpha_3^3 - \alpha_2^2\alpha_3^2 - \alpha_2\alpha_3^3 + 2\alpha_3^4)n^2\}=0. \label{tep7r6}
\end{multline}

Except for the last factor, if we equate any of the other factors on the left-hand side of \eqref{tep7r6}, we get trivial solutions of the diophantine system \eqref{tep7red}. To obtain nontrivial solutions, we equate the last factor to 0, and write 
\begin{equation}
m=ny/(2\alpha_1^2 - \alpha_1\alpha_2 + 5\alpha_1\alpha_3 - \alpha_2^2 - \alpha_2\alpha_3 + 2\alpha_3^2),
\label{valm}
\end{equation}
when we get,
\begin{equation}
y^2=\phi(\alpha_1, \alpha_2, \alpha_3), \label{quartec}
\end{equation}
where 
\begin{equation}
\begin{aligned}
\phi(\alpha_1, \alpha_2, \alpha_3)&=(2\alpha_1^2 - \alpha_1\alpha_2 + 5\alpha_1\alpha_3 - \alpha_2^2 - \alpha_2\alpha_3 + 2\alpha_3^2)
(2\alpha_1^4 - \alpha_1^3\alpha_2 \\ & \quad \;\; + 9\alpha_1^3\alpha_3 - \alpha_1^2\alpha_2^2  + 37\alpha_1^2\alpha_2\alpha_3 + 14\alpha_1^2\alpha_3^2
+ 38\alpha_1\alpha_2^2\alpha_3 \\ & \quad \;\; + 37\alpha_1\alpha_2\alpha_3^2 + 9\alpha_1\alpha_3^3 - \alpha_2^2\alpha_3^2 - \alpha_2\alpha_3^3 + 2\alpha_3^4),\\
&=(\alpha_1^2 - 38\alpha_1\alpha_3 + \alpha_3^2)\alpha_2^4 + 2(\alpha_1 + \alpha_3)(\alpha_1^2 - 38\alpha_1\alpha_3 + \alpha_3^2)\alpha_2^3\\  & \quad \;\; - (3\alpha_1^4 - 26\alpha_1^3\alpha_3 - 98\alpha_1^2\alpha_3^2 - 26\alpha_1\alpha_3^3 + 3\alpha_3^4)\alpha_2^2 \\  & \quad \;\; - 2(2\alpha_1 + \alpha_3)(\alpha_1 + 2\alpha_3)(\alpha_1 + \alpha_3)(\alpha_1^2 - 18\alpha_1\alpha_3 + \alpha_3^2)\alpha_2 \\  & \quad \;\; + (\alpha_1 + 2\alpha_3)^2(2\alpha_1 + \alpha_3)^2(\alpha_1 + \alpha_3)^2.
\end{aligned}
\label{defphi}
\end{equation}

Now $\phi(\alpha_1, \alpha_2, \alpha_3)$ may be considered as a quartic function of $\alpha_2$ in which the constant term is a perfect square, and the obvious method to find solutions of the diophantine equation \eqref{quartec} is to consider it as a quartic model of an elliptic curve over the function field $\mathbb{Q}(\alpha_1, \alpha_3)$. We can readily obtain rational points on this elliptic curve and thus obtain parametric ideal solutions of the TEP of degree 7 but these solutions are of degrees $ \geq 4$.

We will now describe an elementary method of choosing $\alpha_1, \alpha_2, \alpha_3$ such that $\phi(\alpha_1, \alpha_2, \alpha_3)$ becomes a perfect square, and thus obtain parametric ideal solutions of degree 4 of the TEP of degree 7. We first write $\alpha_2=f\alpha_1+g\alpha_3$, where $f$ and $g$ are arbitrary parameters, and consider $\phi(\alpha_i)$ as a homogeneous sextic polynomial in the variables $\alpha_1$ and $\alpha_3$. We choose one of the parameters  (either $f$ or $g$) such that the discriminant of $\phi(\alpha_i)$ with respect to $\alpha_1$ vanishes, and now  $\phi(\alpha_i)$ has a squared factor which may be factored out leaving us with a quartic polynomial in the variables $\alpha_1$ and $\alpha_3$. We now repeat the process with the quartic polynomial, that is, we choose  the remaining parameter (either $f$ or $g$)  such that the discriminant of the quartic polynomial with respect to $
\alpha_1$ vanishes, and now we can again factor out a squared factor from the quartic polynomial leaving us with a quadratic polynomial in the variables $ \alpha_1$ and $\alpha_3$. We can readily equate this quadratic polynomial in  $
\alpha_1$ and $\alpha_3$ to a perfect square, and thus obtain parametric ideal solutions of degree 4 of the TEP of degree 7. 

We will now illustrate the method described above to obtain a specific parametric solution of the TEP of degree 7. On substituting 
$\alpha_2=f\alpha_1+g\alpha_3$ in $\phi(\alpha_1, \alpha_2, \alpha_3)$, and equating to zero  the discriminant, with respect to $\alpha_1$, of the resulting sextic polynomial, we get the following condition:
\begin{multline}
(f + 2)^2(f - 1)^2(f - g)^2(2f - g)^2(f - 2g)^2(f + g + 1)^2(f - 2g - 1)^2\\
\times (2f - g + 1)^2(g + 2)^2(g - 1)^2(9f^2 - 22fg + 9g^2 - 2f - 2g + 9)\\
\times (f^6 + 114f^5g + 4335f^4g^2 + 55100f^3g^3 + 4335f^2g^4 + 114fg^5 + g^6 + 60f^5\\
 + 4620f^4g + 91320f^3g^2 + 91320f^2g^3 + 4620fg^4 + 60g^5 - 20667f^4\\
 + 48228f^3g + 141678f^2g^2 + 48228fg^3 - 20667g^4 - 35620f^3 + 72420f^2g \\
+ 72420fg^2 - 35620g^3 - 20667f^2 + 46254fg - 20667g^2 + 60f + 60g + 1)=0. \label{conddisfirst}
\end{multline}

Eq.~\eqref{conddisfirst} is readily satisfied by choosing $f$ or $g$ such that  one  of the 10 linear factors on the left-hand side of 
Eq.~\eqref{conddisfirst}  becomes zero. As an example, we take $f=-2$, and now $\phi(\alpha_1, \alpha_2, \alpha_3)$ reduces to $\alpha_3^2 \phi_1(\alpha_1, \alpha_3)$, where
\begin{multline}
\phi_1(\alpha_1, \alpha_3)=\{(3g + 7)\alpha_1 - (g + 2)(g - 1)\alpha_3\}\{(3g + 87)\alpha_1^3 - \\
(g^2 + 115g + 64)\alpha_1^2\alpha_3 + (19g + 11)(2g + 1)\alpha_1\alpha_3^2 - (g + 2)(g - 1)\alpha_3^3\}.
\end{multline}

Now on equating to zero the discriminant of $\phi_1(\alpha_1, \alpha_3)$ with respect to $\alpha_1$, we get the following condition:
\begin{multline}
\quad \quad  (g - 1)^4(g + 1)^2(g + 2)^4(g + 3)^2(g + 4)^2(2g + 3)^2\\
\times (g^4 - 226g^3 - 300g^2 - 130g - 155)=0. \quad \quad \quad \quad \label{conddissec}
\end{multline}

This condition is again satisfied quite simply by choosing the value of $g$ in several ways. We take $g=-1$ when $\phi_1(\alpha_1, \alpha_3)$ reduces to $4(2\alpha_1 + \alpha_3)^2(21\alpha_1^2 + 2\alpha_1\alpha_3 + \alpha_3^2)$. We have now reduced our problem to choosing $\alpha_1$ and  $\alpha_3$ such that $21\alpha_1^2 + 2\alpha_1\alpha_3 + \alpha_3^2$ becomes a perfect square. It is readily seen that that this condition is satisfied of we take
\begin{equation}
\alpha_1 = 2t + 2, \alpha_3 = t^2 - 21, \label{vala1a3}
\end{equation}
where $t$ is an arbitrary parameter. Since we took $f=-2$ and $g=-1$, we now get $\alpha_2=-2\alpha_1-g\alpha_3$ and, on using \eqref{vala1a3}, this gives us,
\begin{equation}
\alpha_2 = -t^2 - 4t + 17. \label{vala2}
\end{equation}

With the values of $\alpha_1, \alpha_2, \alpha_3$ given by \eqref{vala1a3} and \eqref{vala2}, the value of $\phi(\alpha_1, \alpha_2, \alpha_3)$ is a perfect square, and hence Eq.~\eqref{tep7r6} can be solved for $m$ and $n$, and finally, using the relations \eqref{RMJMsol}, we get a solution of the simultaneous diophantine equations \eqref{tep7red}. This solution may be written as follows:
\small
\begin{equation}
\begin{aligned}
x_1&=t^4 + 6t^3 - 32t^2 - 158t + 279,\quad & x_2&=4t^3 + 28t^2 + 4t - 420,\\
 x_3&=t^4 + 6t^3 - 4t^2 - 102t - 93,\quad & x_4&=t^4 - 50t^2 - 56t + 393, \\
y_1&=t^4 + 8t^3 - 26t^2 - 112t + 321,\quad & y_2&=t^4 + 2t^3 - 16t^2 + 46t + 63,\\
 y_3&=4t^3 - 4t^2 - 60t + 348,\quad & y_4&=t^4 + 2t^3 - 44t^2 - 10t + 435,
\end{aligned}
\end{equation}
\normalsize
where $t$ is an arbitrary parameter. As a numerical example, on taking $t=2$, we get the solution,
\[
101^r+268^r+249^r+97^r= 73^r+123^r+244^r+271^r, \quad r=2, 4, 6.
\]

We will now find more solutions of the simultaneous diophantine equations \eqref{tep7red} following the method illustrated above. We have already noted that the condition \eqref{conddisfirst} can readily be satisfied in 10 different ways. We consider each possibility,  obtain the quartic polynomial corresponding to the polynomial $\phi_1(\alpha_i)$ of the above illustrative example, and equate its discriminant to zero. Like the condition \eqref{conddissec}, each such condition can be readily satisfied in several ways.  We exhaustively analysed all the  possibilities, and could obtain only four distinct nontrivial solutions of the simultaneous diophantine equations \eqref{tep7red}, one of which has already been given above. 

The remaining three  solutions of degree 4 of the simultaneous diophantine equations \eqref{tep7red} may be written, in terms of an arbitrary parameter $t$,    as follows:
\small
\[
\begin{aligned}
x_1 &= 3t^4 + 40t^3 - 274t^2 + 48t + 1383, &  x_2 &= t^4 - 88t^3 - 214t^2 + 736t - 675, \\
x_3 &= 7t^4 + 14t^3 - 152t^2 + 962t + 609, &  x_4 &= 4t^4 - 70t^3 + 46t^2 + 254t - 1914,\\
y_1 &= 4t^4 - 6t^3 - 274t^2 - 642t + 1158, &  y_2 &= 3t^4 - 92t^3 - 62t^2 + 676t + 1155,\\
y_3 &= 7t^4 + 12t^3 - 14t^2 - 1124t - 81, &  y_4 &= t^4 + 66t^3 + 184t^2 + 238t - 1929,
\end{aligned}
\]
and
\[
\begin{aligned}
x_1 & = 2t^4 + 42t^3 - 170t^2 + 942t - 48, &  x_2 & = t^4 - 84t^3 - 242t^2 - 68t - 1911,\\
 x_3 & = 5t^4 + 12t^3 - 10t^2 + 1724t + 1341, &  x_4 & = 3t^4 - 58t^3 + 92t^2 - 310t - 1263,\\
y_1 & = 3t^4 - 16t^3 - 170t^2 - 1320t - 1569, &  y_2 & = 2t^4 - 86t^3 - 106t^2 - 146t + 1872,\\
 y_3 & = 5t^4 + 10t^3 + 164t^2 - 1562t - 921, &  y_4 & = t^4 + 56t^3 + 266t^2 + 928t - 483, 
\end{aligned}
\]
and
\[
\begin{aligned}
x_1 & = 4t^4 + 18t^3 + 34t^2 + 526t + 378, &  x_2 & = 3t^4 + 16t^3 - 86t^2 - 32t + 579,\\
 x_3 & = t^4 - 12t^3 + 74t^2 + 1020t + 1317, &  x_4 & = 2t^4 + 42t^3 + 110t^2 + 230t + 1056,\\
 y_1 & = 4t^4 + 14t^3 - 122t^2 - 958t - 1338, &  y_2 & =  2t^4 - 14t^3 - 242t^2 - 658t - 48,\\
 y_3 & = 3t^4 + 40t^3 + 74t^2 - 696t - 861, &  y_4 & =t^4 + 44t^3 + 202t^2 + 164t - 891.
\end{aligned}
\]

We have thus obtained four solutions of the simultaneous diophantine equations \eqref{tep7red} in terms of quartic polynomials. Using  the relations \eqref{simplecond}, we immediately get four parametric solutions of degree 4 of the TEP of degree 7.

We note that our solutions, like all other known solutions of the TEP of degree 7, are symmetric solutions, that is, they satisfy the conditions \eqref{simplecond}. It would be of interest to find nonsymmetric solutions, that is, solutions that do not satisfy the simplifying conditions \eqref{simplecond}.

\end{document}